\documentclass[10pt,draft,twoside,reqno,a4paper]{article}
\usepackage{amsmath,amsthm,amsfonts,amssymb,verbatim}
\newtheorem{Theorem}{Theorem}
\newtheorem{Lemma}{Lemma}
\newtheorem{Corollary}{Corollary}
\newtheorem{Proposition}{Proposition}
\newtheorem{Remark}{Remark}
\newcommand{\R}{\mathbb R}
\newcommand{\C}{\mathbb C}
\newcommand{\N}{\mathbb N}

\newcommand{\de}{\partial}
\newcommand{\dbar}{\overline\partial}
\newcommand{\bra}{\langle}
\newcommand{\ket}{\rangle}

\newcommand{\kf}{k_1}
\newcommand{\ks}{k_2}
\newcommand{\tkf}{k_{\tau,1}}
\newcommand{\tks}{k_{\tau,2}}
\newcommand{\tf}{f_\tau}
\newcommand{\mo}{\mu_0}
\newcommand{\no}{\nu_0}
\newcommand{\tmu}{\mu_\tau}
\newcommand{\tnu}{\nu_\tau}
\newcommand{\tmuo}{\mu_{0,\tau}}
\newcommand{\tnuo}{\nu_{0,\tau}}
\newcommand{\iso}{d}
\newcommand{\ef}{\eta_1}
\newcommand{\es}{\eta_2}
\newcommand{\dotef}{\dot\eta_1}
\newcommand{\dotes}{\dot\eta_2}
\renewcommand{\d}{{\mathrm d}}
\begin{document}
\title{On planar Beltrami equations and H\"older regularity}
\author{Tonia Ricciardi\thanks{
Supported in part by 
the MIUR National Project {\em Variational Methods and
Nonlinear Differential Equations}.}\\
{\small Dipartimento di Matematica e Applicazioni
``R.~Caccioppoli"}\\
{\small Universit\`a di Napoli Federico II}\\
{\small Via Cintia, 80126 Napoli, Italy}\\
{\small fax: +39 081 675665}\\
{\small\tt{ tonia.ricciardi@unina.it}}\\
}
\date{}
\maketitle
\begin{abstract}
We provide estimates for the H\"older exponent of solutions to 
the Beltrami equation $\dbar f=\mu\de f+\nu\overline{\de f}$,
where the Beltrami coefficients $\mu,\nu$ satisfy $\||\mu|+|\nu|\|_\infty<1$
and $\Im(\nu)=0$.
Our estimates depend on the arguments of the Beltrami coefficients
as well as on their moduli. 
Furthermore, we 
exhibit a class of mappings of the ``angular stretching" type,
on which our estimates are actually attained, 
and we discuss the main properties of such mappings.
\end{abstract}
\begin{description}
\item {\textsc{Key Words:}} linear Beltrami equation, 
H\"older regularity, angular stretching
\item {\textsc{MSC 2000 Subject Classification:}} 30C62 (35J25)
\end{description}
\section{Introduction and statement of the main\\ results}
\label{sec:intro}
Let $\Omega$ be a bounded open subset of $\R^2$ and let
$f\in W_{\mathrm{loc}}^{1,2}(\Omega,\C)$ satisfy the Beltrami equation
\begin{equation}
\label{beltrami}
\dbar f=\mu\de f+\nu\overline{\de f}
\qquad\mathrm{a.e.\ in\ }\Omega,
\end{equation}
where $\dbar=(\de_1+i\de_2)/2$, $\de=(\de_1-i\de_2)/2$ and
$\mu,\nu$, are bounded, measurable functions  
satisfying $\||\mu|+|\nu|\|_\infty<1$.
Equation \eqref{beltrami} arises in the study of conformal
mappings between domains equipped with measurable Riemannian structures, 
see \cite{AIM}.
By classical work of Morrey~\cite{Mo}, it is well-known that
solutions to \eqref{beltrami} are H\"older continuous.
More precisely, there exists $\alpha\in(0,1)$ such that  
for every compact $T\Subset\Omega$ there exists $C_T>0$
such that
\[
\frac{|f(z)-f(z')|}{|z-z'|^\alpha}\le C_T
\qquad\forall z,z'\in T,\ z\neq z'.
\]
Let 
\[
K_{\mu,\nu}=\frac{1+|\mu|+|\nu|}{1-|\mu|-|\nu|} 
\]
denote the distortion function. 
Then, the following estimate holds:
\begin{equation}
\label{classicalhoelder}
\alpha\ge\|K_{\mu,\nu}\|_\infty^{-1}.
\end{equation}
This estimate is sharp, in the sense that it reduces to
an equality on the radial stretching
\begin{equation}
\label{radialstretch}
f(z)=|z|^{\alpha-1}z,
\end{equation}
which satisfies \eqref{beltrami} with $\mu(z)=-(1-\alpha)/(1+\alpha)z\bar z^{-1}$
and $\nu=0$.
There exists a wide literature concerning the regularity theory
for \eqref{beltrami}, particularly in the degenerate case where
$\||\mu|+|\nu|\|_\infty=1$, or equivalently, when the distortion function 
$K_{\mu,\nu}$ is unbounded. 
See, e.g., \cite{AIKM,GMSV,IM,IS}, and the references therein.
See also \cite{FS}, where an estimate of the constant $C_T$ is given.
Most of the results mentioned above provide estimates 
in terms of the distortion function $K_{\mu,\nu}$,
and there is no loss of generality in assuming that $\nu=0$. Indeed, 
the following ``device of Morrey" may be used, as described in \cite{BN}:
at points where $\overline\de f=0$ we set
$\widetilde\mu=\mu+\nu\overline{\de f}/\dbar f$; at points where $\dbar f=0$ we set
$\widetilde\mu=0$. Then $f$ is a solution to $\dbar f=\widetilde\mu\de f$ and
$|\widetilde\mu|\le|\mu|+|\nu|$.
On the other hand, in this note we are interested in estimates which preserve the
information contained in the arguments of the Beltrami coefficients $\mu,\nu$,
in the spirit of the work of Andreian Cazacu~\cite{Ca} and 
of Reich and Walczak~\cite{RW}.
We restrict our attention to the case $\Im(\nu)=0$. 
This assumption corresponds to assuming that
the Riemannian metric in the target manifold 
is represented by a diagonal matrix-valued function.
We will also show that our estimates are sharp, in the sense that they are
attained in a class of mappings of the ``angular stretching" type
(see ansatz~\eqref{fansatz} below), which generalize the 
radial stretchings \eqref{radialstretch}.
We expect such mappings to be of interest in other 
areas of quasiconformal mapping theory, and therefore we analyze them in some detail. 
Our first result is the following.
\begin{Theorem}
\label{thm:estimate}
Let $f\in W_{\mathrm{loc}}^{1,2}(\Omega,\C)$ satisfy the Beltrami equation~\eqref{beltrami}
with $\Im(\nu)=0$.
Then, $f$ is $\alpha$-H\"older continuous with $\alpha\ge\beta(\mu,\nu)$, where $\beta(\mu,\nu)$
is defined by
\begin{align}
\label{betaestimate}
\beta&(\mu,\nu)^{-1}=
\sup_{S_\rho(x)\subset\Omega}
\inf_{\varphi,\psi\in\mathcal B_{x,\rho}}\sqrt{\frac{\sup\varphi}{\inf\psi}}\\
\nonumber
&\Big\{\frac{1}{|S_\rho(x)|}\int_{S_\rho(x)}\sqrt{\frac{\psi}{\varphi}}
\frac{|1-\overline n^2\mu|^2-\nu^2}
{\sqrt{1-(|\mu|+\nu)^2}\sqrt{1-(|\mu|-\nu)^2}}\d\sigma\times\\
\nonumber
&\qquad\times\left(\frac{4}{\pi}\arctan\left(\frac{\inf_{S_\rho(x)}
\frac{(1-\nu)^2-|\mu|^2}{(1+\nu)^2-|\mu|^2}/\varphi\psi}
{\sup_{S_\rho(x)}\frac{(1-\nu)^2-|\mu|^2}{(1+\nu)^2-|\mu|^2}/\varphi\psi}\right)^{1/4}\right)^{-1}\Big\}.
\end{align}
Here $\mathcal B_{x,\rho}$ denotes the set of positive functions
in $L^\infty(S_\rho(x))$ which are bounded below away from zero, and
$n$ denotes complex number corresponding to the outer unit normal
to $S_\rho(x)$. 
\end{Theorem}
Estimate \eqref{betaestimate} improves the classical estimate~\eqref{classicalhoelder};
a verification is provided in Section~\ref{sec:proofs}, Remark~\ref{rem:check}.
In Theorem~\ref{thm:sharp} below 
we will show that estimate~\eqref{betaestimate} is sharp, in the sense that it reduces
to an equality when $\mu,\nu$
are of the special form
\[
\mu(z)=-\mu_0(\arg z)z\overline z^{-1},
\qquad \nu(z)=-\nu_0(\arg z)
\]
and $f$ is of the ``angular stretching" form
\[
f(z)=|z|^\alpha(\eta_1(\arg z)+i\eta_2(\arg z)),
\]
for suitable choices of the bounded, $2\pi$-periodic functions 
$\mu_0,\nu_0,\eta_1,\eta_2:\R\to\R$.
The following weaker form of estimate~\eqref{betaestimate}
is obtained by taking
$\varphi=\psi=1$. 
\begin{Corollary}
\label{cor:phipsi=1}
Let $f\in W_{\mathrm{loc}}^{1,2}(\Omega,\C)$ satisfy the Beltrami equation~\eqref{beltrami}
with $\Im(\nu)=0$.
Then, $f$ is $\alpha$-H\"older continuous with
\begin{align}
\label{alphaestimate}
\alpha\ge
\left\{\sup_{S_\rho(x)\subset\Omega}
\frac{\frac{1}{|S_\rho(x)|}\int_{S_\rho(x)}
\frac{|1-\overline n^2\mu|^2-\nu^2}
{\sqrt{1-(|\mu|+\nu)^2}\sqrt{1-(|\mu|-\nu)^2}}\d\sigma}
{\frac{4}{\pi}\arctan\left(\frac{\inf_{S_\rho(x)}
\frac{(1-\nu)^2-|\mu|^2}{(1+\nu)^2-|\mu|^2}}
{\sup_{S_\rho(x)}\frac{(1-\nu)^2-|\mu|^2}{(1+\nu)^2-|\mu|^2}}\right)^{1/4}}\right\}^{-1}.
\end{align}
\end{Corollary}
This estimate is also sharp, in the sense that it actually
reduces to an equality for suitable choices of $\mu,\nu$ and $f$, 
but it does not contain estimate~\eqref{classicalhoelder} as a special case. 
We now show that estimate~\eqref{alphaestimate} contains 
some known results for $\mu=0$ and for $\nu=0$
as special cases.
\paragraph{Special case $\nu=0$.}
This case corresponds to assuming that the target domain
is equipped with the standard Euclidean metric.
In this special case, our estimate yields
\begin{equation}
\label{specialalphaestimate}
\alpha\ge\left\{\sup_{S_\rho(x)\subset\Omega}
\frac{1}{|S_\rho(x)|}\int_{S_\rho(x)}
\frac{|1-\overline n^2\mu |^2}{1-|\mu|^2}\,\mathrm{d}\sigma\right\}^{-1},
\end{equation}
which may also be obtained from the estimate in \cite{sharpholder} for elliptic equations
whose coefficient matrix has unit determinant.
We note that the integrand function
\begin{align*}
\frac{|1-\overline n^2\mu |^2}{1-|\mu|^2}
=\frac{|D_{\overline n}f|^2}{J_f}=
K_{\mu,0}-2\frac{|\mu|+\Re\left(\mu,n^2\right)}{1-|\mu|^2}
\end{align*}
also appears in \cite{RW} in the study of the conformal 
modulus of rings.
\paragraph{Special case $\mu=0$.}
This case corresponds to assuming that
the metric on $\Omega$ is Euclidean. In this case,
estimate \eqref{alphaestimate} yields 
\begin{align}
\label{thirdalphaestimate}
\alpha\ge\sup_{S_\rho(x)\subset\Omega}
\frac{4}{\pi}\arctan\left(\frac{\inf_{S_\rho(x)}\frac{1-\nu}{1+\nu}}
{\sup_{S_\rho(x)}\frac{1-\nu}{1+\nu}}\right)^{1/2}
\ge\frac{4}{\pi}\arctan\|K\|_\infty^{-1},
\end{align}
which is a consequence of the sharp H\"older estimate obtained in
Piccinini and Spagnolo~\cite{PS} for isotropic elliptic equations.
\par
In Theorem~\ref{thm:sharp} below we assert that
the equality $\alpha=\beta(\mu,\nu)$ 
may hold even when \textit{both} $\mu\neq0$ and $\nu\neq0$. 
We denote by $B$ the unit disk in $\R^2$.
\begin{Theorem}
\label{thm:sharp}
For every $\tau\in[0,1]$ there exist $\alpha_\tau>0$, $2\pi$-periodic functions 
$\eta_{\tau,1},\eta_{\tau,2}\in W_{\mathrm{loc}}^{1,2}(\R)$
and corresponding coefficients $\tmu,\tnu$, depending on the angular variable only,
with the following properties:
\begin{enumerate}
\item[(i)]
The mapping $\tf\in W_{\mathrm{loc}}^{1,2}(B)$ defined in $B\setminus\{0\}$ by
\begin{equation*}
\tf(z)=|z|^{\alpha_\tau}\left(\eta_{\tau,1}(\arg z)+i\eta_{\tau,2}(\arg z)\right)
\end{equation*}
satisfies \eqref{beltrami} with $\mu=\tmu$ and $\nu=\tnu$;
\item[(ii)]
$\beta(\tmu,\tnu)=\alpha_\tau$;
\item[(iii)]
$\tmu=0$ if and only if $\tau=0$; $\tnu=0$
if and only if $\tau=1$.
\end{enumerate}
\end{Theorem}
This note is organized as follows.
In Section~\ref{sec:prelims} we derive the basic properties of
the mappings of the ``angular stretching" form,
which naturally appear in our problem. 
In Section~\ref{sec:proofs}
we provide the proofs of Theorem~\ref{thm:estimate}  and Theorem~\ref{thm:sharp}, 
which are based on the equivalence
between Beltrami equations and elliptic divergence form equations, 
to which we can apply some recent results in \cite{besthoelder}.
A proof of the equivalence is provided in the Appendix.
\section{Angular stretchings}
\label{sec:prelims}
In this section we derive some properties of functions of the
``angular stretching" form:
\begin{equation}
\label{fansatz}
f(z)=|z|^\alpha(\eta_1(\arg z)+i\eta_2(\arg z)),
\end{equation}
where $\alpha\in\R$ and $\eta_1,\eta_2:\R\to\R$ are $2\pi$-periodic functions.
We assume $\alpha>0$ and $\ef,\es\in W_{\mathrm{loc}}^{1,2}(\R)$ so that
$f\in W_{\mathrm{loc}}^{1,2}(\C)$.
We note that mappings of the form \eqref{fansatz} generalize the radial stretchings
\eqref{radialstretch}.
We also note that $f$ is injective if and only if $\ef^2(\theta)+\es^2(\theta)\neq0$
for all $\theta\in\R$, $\ef,\es$ have minimal period $2\pi$ 
and $\dotes\ef-\dotef\es=(\ef^2+\es^2)(\mathrm{d}/\mathrm{d}\theta)\arg(\ef+i\es)$ 
has constant sign a.e. 
Recalling that in polar coordinates $x=r\cos\theta$, $y=r\sin\theta$ we have
\[
\left[\begin{matrix}\de_x\\
\de_y\end{matrix}\right]
=\left[\begin{matrix}
\cos\theta&-\sin\theta\\
\sin\theta&\cos\theta
\end{matrix}\right]
\left[\begin{matrix}\de_r\\
\frac{1}{r}\de_\theta\end{matrix}\right],
\]
we derive, at every point in $\R^2\setminus\{0\}$:
\begin{equation*}
r^{-(\alpha-1)}Df=
\left[\begin{matrix}
\alpha\cos\theta\ef-\sin\theta\dot\ef
&&\alpha\sin\theta\ef+\cos\theta\dot\ef\\
\alpha\cos\theta\es-\sin\theta\dot\es
&&\alpha\sin\theta\es+\cos\theta\dot\es
\end{matrix}\right]
\end{equation*}
so that the Jacobian $J_f$ is given by
\begin{equation}
\label{Jf}
r^{-2(\alpha-1)}J_f=\alpha(\ef\dot\es-\dot\ef\es).
\end{equation}
Since for any $2\times2$ matrix $A=(a_{ij})$, $j=1,2$,
we have $\det(AA^T)=\det(A^TA)=(\det A)^2$
and $\mathrm{tr}\,(AA^T)=\mathrm{tr}\,(A^TA)=\sum_{i,j=1}^2a_{ij}^2$,
the operator norm of $Df$ equals the operator norm of $Df^T$, which is more easily calculated.
The tensor $DfDf^T$ is given by
\begin{equation*}
\label{DfDf}
r^{-2(\alpha-1)}DfDf^T=
\left[\begin{matrix}
\alpha^2\ef^2+\dot\ef^2&&\alpha^2\ef\es+\dot\ef\dot\es\\
\alpha^2\ef\es+\dot\ef\dot\es&&\alpha^2\es^2+\dot\es^2
\end{matrix}\right].
\end{equation*}
Hence,
\begin{align*}
\mathrm{tr}\,(r^{-2(\alpha-1)}DfDf^T)=&\alpha^2(\eta_1^2+\eta_2^2)+\dotef^2+\dotes^2\\
\det(r^{-2(\alpha-1)}DfDf^T)=&\alpha^2(\ef\dotes-\dotef\es)^2
\end{align*}
and therefore the eigenvalues of $r^{-2(\alpha-1)}Df^TDf$ are the solutions to
the second order algebraic equation
\begin{equation*}
\lambda^2-[\alpha^2(\ef^2+\es^2)+\dot\ef^2+\dot\es^2]\lambda
+\alpha^2(\ef\dot\es-\dot\ef\es)^2=0.
\end{equation*}
The discriminant of the equation above is given by
\[
\mathcal D=[\alpha^2(\ef^2+\es^2)+\dot\ef^2+\dot\es^2]^2-4\alpha^2(\ef\dot\es-\dot\ef\es)^2.
\]
In view of the elementary identity 
$(a^2+b^2+c^2+d^2)^2-4(ad-bc)^2=(a^2-d^2)^2+(b^2-c^2)^2+2(ab+cd)^2+2(ac+bd)^2$
for every $a,b,c,d\in\R$, we derive the equivalent expression
\begin{align*}
\label{discriminant}
\mathcal D=(\alpha^2\ef^2-\dotes^2)^2+&(\alpha^2\es^2-\dotef^2)^2
+2(\alpha^2\ef\es+\dotef\dotes)^2\\+&2\alpha^2(\ef\dotef+\es\dotes)^2.
\end{align*}
We may write
\begin{align*}
r^{-2(\alpha-1)}|Df|^2=
\frac{1}{2}\left\{\alpha^2(\ef^2+\es^2)+\dot\ef^2+\dot\es^2+\sqrt{\mathcal D}\right\}.
\end{align*}
Therefore, at every point in $\R^2\setminus\{0\}$ the distortion of $f$ is given by
\begin{align*}
\frac{|Df|^2}{J_f}=\frac{\alpha^2(\ef^2+\es^2)+\dot\ef^2+\dot\es^2+\sqrt{\mathcal D}}
{2\alpha(\ef\dot\es-\dot\ef\es)}.
\end{align*}
In particular, $f$ has bounded distortion if and only if
\begin{equation*}
\ef^2+\es^2+\dotef^2+\dotes^2\le C(\ef\dot\es-\dot\ef\es)
\end{equation*}
for some constant $C>0$.
\par
In order to prove Theorem~\ref{thm:sharp} we need some properties for the special case
where $f$ is of the angular stretching form \eqref{fansatz} 
and moreover $f$ satisfies the Beltrami equation~\eqref{beltrami} 
with $\mu,\nu$ of the special form
\begin{equation}
\label{mu}
\mu(z)=-\mu_0(\arg z)\,z\bar z^{-1}
\end{equation}
and
\begin{equation}
\label{nu}
\nu(z)=-\nu_0(\arg z),
\end{equation}
for some bounded, $2\pi$-periodic functions 
$\mo,\no$ such that $\||\mo|+|\no|\|_\infty<1$.
We use the following facts.
\begin{Proposition}
\label{prop:polar}
Suppose $f$ is of the angular stretching form~\eqref{fansatz}
and satisfies the Beltrami equation~\eqref{beltrami} 
with $\mu,\nu$ given by \eqref{mu}--\eqref{nu}.
Then, $(\eta_1,\eta_2)$ 
satisfies the system:
\begin{align}
\label{system}
\begin{cases}
\dotef=-\alpha\ks^{-1}\es\\
\dotes=\alpha\kf\ef,
\end{cases}
\end{align}
where $\kf,\ks>0$ are defined by
\begin{align}
\label{k}
\kf=\frac{1+\mo+\no}{1-\mo-\no},
&&\ks=\frac{1-\mo+\no}{1+\mo-\no}.
\end{align}
Conversely, if $(\ef,\es)$ satisfies \eqref{system} for some
$\alpha>0$ and for some $2\pi$-periodic
functions $\kf,\ks>0$ bounded from above and from below away from zero, 
then $f$ defined by \eqref{fansatz} is a solution to
\eqref{beltrami} with $\mu,\nu$ defined in \eqref{mu}--\eqref{nu}
and $\mo,\no$ given by
\begin{align}
\label{mono}
\mo=\frac{\kf-\ks}{1+\kf+\ks+\kf\ks},
&&\no=\frac{\kf\ks-1}{1+\kf+\ks+\kf\ks}.
\end{align}
\end{Proposition}
\begin{proof}
In polar cooordinates
$x=r\cos\theta$, $y=r\sin\theta$ we have 
\begin{align*}
\dbar=&\frac{1}{2}(\de_x+i\de_y)=\frac{e^{i\theta}}{2}\left(\de_r
+i\,\frac{\de_\theta}{r}\right)\\
\de=&\frac{1}{2}(\de_x-i\de_y)=\frac{e^{-i\theta}}{2}\left(\de_r
-i\,\frac{\de_\theta}{r}\right).
\end{align*}
Hence, \eqref{beltrami} is equivalent to
\begin{equation*}
(e^{i\theta}-\mu e^{-i\theta})f_r-\nu e^{i\theta}\overline{f_r}=
-\frac{i}{r}\left[(e^{i\theta}+\mu e^{-i\theta})f_\theta
-\nu e^{i\theta}\overline{f_\theta}\right].
\end{equation*}
In view of the form~\eqref{mu} of $\mu$
and the of form~\eqref{nu} of $\nu$, the equation above
is equivalent to
\begin{equation*}
(1+\mo)f_r+\no\overline{f_r}=-\frac{i}{r}[(1-\mo)f_\theta+\no\overline{f_\theta}].
\end{equation*}
We compute
\begin{align*}
f_r=\alpha r^{\alpha-1}(\eta_1+i\eta_2),
&&f_\theta=r^\alpha(\dot\eta_1+i\dot\eta_2).
\end{align*}
Substitution yields
\begin{equation}
\label{fpolar}
\alpha(1+\mo+\no)\eta_1+i\alpha(1+\mo-\no)\eta_2
=(1-\mo-\no)\dot\eta_2-i(1-\mo+\no)\dot\eta_1.
\end{equation}
Hence, $(\eta_1,\eta_2)$ satisfies the system~\eqref{system},
with $\kf,\ks$ defined by \eqref{k}.
Conversely, suppose $(\eta_1,\eta_2)$ satisfies \eqref{system}
for some $2\pi$-periodic
functions $\kf,\ks>0$ bounded from above and from below away from zero 
and for some $\alpha>0$.
Then the functions $\mo,\no$ such that \eqref{k} is satisfied
are uniquely defined by \eqref{mono} as the solutions to the linear system
\begin{align*}
(1+\kf)\mo+(1+\kf)\no=&-1+\kf\\
-(1+\ks)\mo+(1+\ks)\no=&-1+\ks.
\end{align*} 
It follows that \eqref{system} is equivalent to \eqref{fpolar},
with $f$ defined by \eqref{fansatz}.
\end{proof}
We finally observe that if $(\ef,\es)$ is
a solution of the system~\eqref{system}, then
the Jacobian determinant of $f$ is given by
\begin{align*}
r^{-2(\alpha-1)}J_f=&\alpha^2(\kf\ef^2+\ks^{-1}\es^2)
\end{align*}
and furthermore, 
\begin{align}
\label{knorm}
\frac{|Df|^2}{J_f}
=&\Big[2(\kf\ef^2+\ks^{-1}\es^2)\Big]^{-1}\Big[(1+\kf^2)\ef^2+(1+\ks^{-2})\es^2+\\
\nonumber
&+\sqrt{(1-\kf^2)^2\ef^4+(1-\ks^{-2})^2\es^4+2[(1-\kf\ks^{-1})^2
+(\kf-\ks^{-1})^2]\ef^2\es^2}\Big].
\end{align}
We also note that system~\eqref{system} implies that $\ef$ is a $2\pi$-periodic solution to
the weighted Sturm-Liouville equation
\[
\frac{\mathrm{d}}{\mathrm{d}t}(\ks\dotef)+\alpha^2\kf\ef=0
\]
and similarly $\es$ satisfies
\[
\frac{\mathrm{d}}{\mathrm{d}t}(\kf^{-1}\dotes)+\alpha^2\ks^{-1}\es=0.
\]
\subparagraph{Special case $\nu=0$.}
The results described in Proposition~\ref{prop:polar} take a particularly simple form when $\nu=0$,
which is equivalent to $\kf=\ks^{-1}=k$. In this case system~\eqref{system}
reduces to
\begin{equation}
\label{specialsys}
\begin{cases}
\dotef=-\alpha k\es\\
\dotes=\alpha k\ef
\end{cases}
\end{equation}
which may be explicitly solved. Indeed, from
\eqref{specialsys} 
we derive $\dotef\ef+\dotes\es=0$ and therefore $\ef^2+\es^2$ is constant.
By linearity we may assume $\ef^2+\es^2\equiv1$. Hence, there exists a funtion $h(\theta)$
such that $\ef(\theta)=\cos h(\theta)$ and $\es(\theta)=\sin h(\theta)$. By \eqref{specialsys}
we conclude that up to an additive constant $h(\theta)=\alpha\int_0^\theta k$,
and therefore we obtain that $f(z)=|z|^\alpha\exp\{i\alpha\int_0^\theta k\}$.
In view of the $2\pi$-periodicity of $\ef,\es$ we also obtain that
$\alpha=2\pi n(\int_0^{2\pi}k)^{-1}$ for some $n\in\N$.
From equation~\eqref{knorm} we derive, for every $z\neq0$:
\[
\frac{|Df|^2}{J_f}=\frac{1+k^2+|1-k^2|}{2k}=\max\{k,k^{-1}\}.
\]
Since $k\ge1$ if and only if $\mo\ge0$, the expression above implies 
the known fact
\[
\frac{|Df|^2}{J_f}=\frac{1+|\mu|}{1-|\mu|}=K_{\mu,0}.
\]
\section{Proofs }
\label{sec:proofs}
We first of all check that estimate~\eqref{betaestimate} 
in Theorem~\ref{thm:estimate} actually improves
the classical estimate~\eqref{classicalhoelder}.
\begin{Remark}
\label{rem:check}
The following estimate holds:
\begin{align*}
\beta(\mu,\nu)\ge\|K_{\mu,\nu}\|_\infty^{-1},
\end{align*}
where $\beta(\mu,\nu)$ is the quantity defined in Theorem~\ref{thm:estimate}.
\end{Remark}
\begin{proof}
Recall from Section~\ref{sec:intro} that
$K_{\mu,\nu}=(1+|\mu|+|\nu|)/(1-|\mu|-|\nu|)$.
For every $S_\rho(x)\subset\Omega$, we choose
\begin{align*}
&\varphi=\left.\frac{|1-\overline n^2\mu|^2-\nu^2}{(1+\nu)^2-|\mu|^2}\,\right|_{S_\rho(x)},
&&\psi=\left.\frac{(1-\nu)^2-|\mu|^2}{|1-\overline n^2\mu|^2-\nu^2}\,\right|_{S_\rho(x)}.
\end{align*}
We have that
\begin{align*}
\sup\varphi\le&\sup\frac{(1+|\mu|)^2-\nu^2}{(1+\nu)^2-|\mu|^2}
=\sup\frac{1+|\mu|-\nu}{1-|\mu|+\nu}\le\|K_{\mu,\nu}\|_\infty\\
\inf\psi\ge&\inf\frac{(1-\nu)^2-|\mu|^2}{(1+|\mu|)^2-\nu^2}
=\inf\frac{1-|\mu|-\nu}{1+|\mu|+\nu}\ge\|K_{\mu,\nu}\|_\infty^{-1}
\end{align*}
and therefore
\[
\frac{\sup\varphi}{\inf\psi}\le\|K_{\mu,\nu}\|_\infty^2.
\]
Moreover,
\begin{align*}
\varphi\psi=&\left.\frac{(1-\nu)^2-|\mu|^2}{(1+\nu)^2-|\mu|^2}\right|_{S_\rho(x)}.
\end{align*}
In view of the elementary identity
\[
[(1-\nu)^2-|\mu|^2][(1+\nu)^2-|\mu|^2]
=[1-(|\mu|+\nu)^2][1-(|\mu|-\nu)^2]
\]
we finally obtain
\begin{align*}
\frac{\psi}{\varphi}=\left.\frac{(1-(|\mu|+\nu)^2)(1-(|\mu|-\nu)^2)}
{(|1-\overline n^2\mu|^2-\nu^2)^2}\right|_{S_\rho(x)}.
\end{align*}
Consequently, inserting into \eqref{betaestimate}, we find
that for every $S_\rho(x)\subset\Omega$:
\begin{align*}
\inf_{\varphi,\psi\in\mathcal B_{x,\rho}}&\sqrt{\frac{\sup\varphi}{\inf\psi}}
\Big\{\frac{1}{|S_\rho(x)|}\int_{S_\rho(x)}\sqrt{\frac{\psi}{\varphi}}
\frac{|1-\overline n^2\mu|^2-\nu^2}
{\sqrt{1-(|\mu|+\nu)^2}\sqrt{1-(|\mu|-\nu)^2}}\d\sigma\times\\
\nonumber
&\times\left(\frac{4}{\pi}\arctan\left(\frac{\inf_{S_\rho(x)}
\frac{(1-\nu)^2-|\mu|^2}{(1+\nu)^2-|\mu|^2}/\varphi\psi}
{\sup_{S_\rho(x)}\frac{(1-\nu)^2-|\mu|^2}{(1+\nu)^2-|\mu|^2}/\varphi\psi}\right)^{1/4}\right)^{-1}\Big\}
\le\|K_{\mu,\nu}\|_{\infty}.
\end{align*}
Consequently,
\[
\beta(\mu,\nu)^{-1}\le\|K_{\mu,\nu}\|_{\infty},
\]
and the asserted estimate is verified.
\end{proof}
We use some results in
\cite{besthoelder} for solutions to the elliptic divergence form equation
\begin{equation}
\label{elliptic}
\mathrm{div}(A\nabla\cdot)=0\qquad\mathrm{in\ }\Omega
\end{equation}
where $A$ is a bounded and symmetric matrix-valued function.
More precisely, let 
\begin{equation}
\label{J}
J(\theta)=\left[\begin{matrix}
\cos\theta&-\sin\theta\\
\sin\theta&\cos\theta
\end{matrix}\right].
\end{equation}
For every $M>1$, let 
\begin{align}
\label{c}
c=c(M,\tau)=\frac{2}{1+M^{-\tau}},
&&\iso=\iso(M,\tau)=\frac{4}{\pi}\arctan M^{-(1-\tau)/2}.
\end{align}
Note that when $\tau=0$ we have $\iso=4\pi^{-1}\arctan M^{-1/2}$
and $c=1$, and when $\tau=1$ we have $\iso=1$ and $c=2/(1+M^{-1})$.
We define the intervals
\[
I_1=[0,\frac{c\pi}{2}),\qquad
I_2=[\frac{c\pi}{2},\pi),\qquad
I_3=[\pi,\pi+\frac{c\pi}{2}),\qquad
I_4=[\pi+\frac{c\pi}{2},2\pi).
\]
Let $\Theta_{\tau,1},\Theta_{\tau,2}:\R\to\R$ be the $2\pi$-periodic Lipschitz functions defined 
in $[0,2\pi)$ by
\[
\Theta_{\tau,1}(\theta)=\begin{cases}
\sin[\iso(c^{-1}\theta-\pi/4)],&\theta\in I_1\\
M^{-(1-\tau)/2}\cos[\iso(c^{-1}M^\tau(\theta-c\pi/2)-\pi/4)],&\theta\in I_2\\
-\sin[\iso(c^{-1}(\theta-\pi)-\pi/4)],&\theta\in I_3\\
-M^{-(1-\tau)/2}\cos[\iso(c^{-1}M^\tau(\theta-\pi-c\pi/2)-\pi/4)],&\theta\in I_4
\end{cases}
\]
and
\[
\Theta_{\tau,2}(\theta)=\begin{cases}
-\cos[\iso(c^{-1}\theta-\pi/4)],&\theta\in I_1\\
M^{(1-\tau)/2}\sin[\iso(c^{-1}M^\tau(\theta-c\pi/2)-\pi/4)],&\theta\in I_2\\
\cos[\iso(c^{-1}(\theta-\pi)-\pi/4)],&\theta\in I_3\\
-M^{(1-\tau)/2}\sin[\iso(c^{-1}M^\tau(\theta-\pi-c\pi/2)-\pi/4)],&\theta\in I_4.
\end{cases}
\]
The following facts were established in \cite{besthoelder}
and will be used in the sequel.
\begin{Theorem}[\cite{besthoelder}]
\label{thm:besthoelder}
The following estimates hold.
\begin{enumerate}
\item[(i)]
Let $w\in W_{\mathrm{loc}}^{1,2}(\Omega)$ be a weak solution to \eqref{elliptic}.
Then, $w$ is $\alpha$-H\"older continuous with $\alpha\ge\gamma(A)$, where
\begin{equation}
\label{gamma}
\gamma(A)=\left(
\sup_{S_\rho(x)\subset\Omega}\inf_{\varphi,\psi\in\mathcal B_{x,\rho}}
\sqrt{\frac{\sup\varphi}{\inf\psi}}\frac{\frac{1}{|S_\rho(x)|}\int_{S_\rho(x)}
\sqrt{\frac{\psi}{\varphi}}\frac{\bra n,A\,n\ket}{\sqrt{\det A}}}
{\frac{4}{\pi}\arctan\left(\frac{\inf_{S_\rho(x)}\det A/\varphi\psi}
{\sup_{S_\rho(x)}\det A/\varphi\psi}\right)^{1/4}}\right)^{-1}
\end{equation}
and where $n$ denotes the outer unit normal.
\item[(ii)]
For every $\tau\in[0,1]$ let $A_\tau$ be the symmetric matrix-valued function 
defined for every $z\neq0$ by
\begin{equation}
\label{tildeA}
A_\tau(z)=
(\tkf(\arg z)-\tks(\arg z))\frac{z\otimes z}{|z|^2}+\tks(\arg z)\mathbf I,
\end{equation}
where $\tkf,\tks$ piecewise constant, $2\pi$-periodic functions defined by
\begin{equation}
\label{tkf}
\tkf(\theta)=\begin{cases}
1,&\mathrm{if\ }\theta\in I_1\cup I_3\\
M,&\mathrm{if\ }\theta\in I_2\cup I_4
\end{cases}
\end{equation}
and
\begin{equation}
\label{tks}
\tks(\theta)=\begin{cases}
1,&\mathrm{if\ }\theta\in I_1\cup I_3\\
M^{1-2\tau},&\mathrm{if\ }\theta\in I_2\cup I_4.
\end{cases}
\end{equation}
There exists $M_0>1$ such that
\begin{equation}
\label{supso}
\gamma(A_\tau)=\frac{\iso}{c},
\end{equation}
for every $M\in(1,M_0^{1/\tau})$, if $\tau>0$,
and with no restriction on $M$ if $\tau=0$.
Furthermore, the function $u_\tau=|z|^{\iso/c}\,\Theta_1(\arg z)$
is a weak solution to \eqref{elliptic} with $A=A_\tau$.
\end{enumerate}
\end{Theorem}
We note that the matrix $A_\tau$ may be equivalently written in the form
\begin{align*}
A_\tau(z)=&\left[\begin{matrix}
\tkf\cos^2\theta+\tks\sin^2\theta&(\tkf-\tks)\sin\theta\cos\theta\\
(\tkf-\tks)\sin\theta\cos\theta&\tkf\sin^2\theta+\tks\cos^2\theta
\end{matrix}\right]\\
=&JK_\tau J^T
\end{align*}
where $K_\tau=\mathrm{diag}(\tkf,\tks)$.
The equivalence between 
Beltrami equations and
elliptic equations 
of the form~\eqref{elliptic}
is well-known.
Indeed, for every matrix $A$ let
\begin{equation}
\label{hatA}
\widehat A=\frac{A}{\det A}.
\end{equation}
The following result holds.
\begin{Lemma}
\label{lem:reduction}
Let $f\in W_{\mathrm{loc}}^{1,2}(\Omega,\C)$ be a solution to \eqref{beltrami} 
with $\Im(\nu)=0$ and
let $A_{\mu,\nu}$ be defined by
\begin{align}
\label{Amu}
A_{\mu,\nu}=\frac{1}{\Delta}\left(\left[\begin{matrix}
|1-\mu|^2&-2\Im(\mu)\\
-2\Im(\mu)&|1+\mu|^2
\end{matrix}\right]-\nu^2\mathbf{I}\right),
\end{align}
where $\Delta=(1+|\mu|+\nu)(1-|\mu|+\nu)$.
Then,
$\Re(f)$ satisfies \eqref{elliptic} with $A=A_{\mu,\nu}$ and
$\Im(f)$ satisfies \eqref{elliptic} with $A=\widehat A_{\mu,\nu}$.
\end{Lemma}
A proof of 
Lemma~\ref{lem:reduction} is provided in the Appendix.
\begin{Lemma}
\label{lem:hat}
For any matrix valued function $A$ we have
\[
\gamma(A)=\gamma(\widehat A)
\]
where $\gamma(A)$ is the quantity defined in Theorem~\ref{thm:besthoelder}.
\end{Lemma}
\begin{proof}
We have $\det\widehat A=(\det A)^{-1}$, and therefore
\begin{equation}
\label{fracA}
\frac{\widehat A}{\sqrt{\det\widehat A}}
=\frac{A}{\sqrt{\det A}}.
\end{equation}
Furthermore, for every $S\subset\Omega$
and for every $\varphi,\psi\in L^\infty(S)$,
\[
\frac{\sup\varphi}{\inf\psi}=\frac{\sup\psi^{-1}}{\inf\varphi^{-1}}
\]
and
\begin{align*}
&\inf_S\frac{\det\widehat A}{\varphi\psi}=\frac{1}{\sup_S(\varphi\psi\det A)},
&&\sup_S\frac{\det\widehat A}{\varphi\psi}=\frac{1}{\inf_S(\varphi\psi\det A)}.
\end{align*}
Hence,
\begin{equation}
\label{infsup}
\frac{\inf_S\det\widehat A/(\varphi\psi)}{\sup_S\det\widehat A/(\varphi\psi)}
=\frac{\inf_S\det A/(\varphi^{-1}\psi^{-1})}{\sup_S\det A/(\varphi^{-1}\psi^{-1})}.
\end{equation}
It follows from \eqref{fracA} and \eqref{infsup} that for any function $F:\R\to\R$
\begin{align*}
&\sqrt{\frac{\sup\varphi}{\inf\psi}}\frac{1}{|S_{\rho}(x)|}\int_{S_\rho(x)}
\sqrt{\frac{\psi}{\varphi}}\frac{\bra n,\widehat An\ket}{\sqrt{\det\widehat A}}\
F\left(\frac{\inf_{S_\rho(x)}\frac{\det\widehat A}{\varphi\psi}}
{\sup_{S_\rho(x)}\frac{\det\widehat A}{\varphi\psi}}\right)\\
=&\sqrt{\frac{\sup\psi^{-1}}{\inf\varphi^{-1}}}\frac{1}{|S_{\rho}(x)|}\int_{S_\rho(x)}
\sqrt{\frac{\varphi^{-1}}{\psi^{-1}}}\frac{\bra n,An\ket}{\sqrt{\det A}}\
F\left(\frac{\inf_{S_\rho(x)}\frac{\det A}{\varphi^{-1}\psi^{-1}}}
{\sup_{S_\rho(x)}\frac{\det A}{\varphi^{-1}\psi^{-1}}}\right).
\end{align*}
Now the statement follows by taking $F(t)=(4\pi^{-1}\arctan t^{1/4})^{-1}$
and observing that $\varphi^{-1}\in\mathcal B_{x,\rho}$ whenever
$\varphi\in\mathcal B_{x,\rho}$.
\end{proof}
\begin{proof}[Proof of Theorem~\ref{thm:estimate}]
In view of Lemma~\ref{lem:reduction}, Lemma~\ref{lem:hat} and Theorem~\ref{thm:besthoelder},
$\Re(g)$ and $\Im(g)$ are $\alpha$-H\"older continuous with
$\alpha\ge\gamma(A_{\mu,\nu})$,
where $A_{\mu,\nu}$ is the matrix defined in \eqref{Amu}.
Setting $\xi=x+\rho e^{it}$, $t\in\R$ for every $\xi\in S_\rho(x)\subset\Omega$, 
we have $n(\xi)=e^{it}$.
We recall that $\Delta=(1+|\mu|+\nu)(1-|\mu|+\nu)=(1+\nu)^2-|\mu|^2$.
Hence, we compute
\begin{align*}
\Delta\,\bra n(\xi),&A_{\mu,\nu}(\xi)n(\xi)\ket=\Delta\bra e^{it},A_{\mu,\nu}(\xi)e^{it}\ket\\
=&\Delta\,(a_{11}\cos^2t+2a_{12}\sin t\cos t+a_{22}\sin^2t)\\
=&1+|\mu|^2-\nu^2-2(\Re(\mu)\cos2t+\Im(\mu)\sin2t)=|1-\overline n^2\mu|^2-\nu^2.
\end{align*}
Furthermore,
\begin{align*}
\Delta^2\det&A_{\mu,\nu}
=(|1-\mu|^2-\nu^2)(|1+\mu|^2-\nu^2)-4\Im(\mu)^2\\
=&(1+|\mu|^2-\nu^2)^2-4|\mu|^2
=((1-|\mu|)^2-\nu^2)((1+|\mu|)^2-\nu^2)\\
=&(1-|\mu|+\nu)(1-|\mu|-\nu)(1+|\mu|+\nu)(1+|\mu|-\nu)\\
=&(1-(|\mu|-\nu)^2)(1-(|\mu|+\nu)^2)
\end{align*}
and therefore
\[
\frac{\bra n,A_{\mu,\nu}n\ket}{\sqrt{\det A_{\mu,\nu}}}
=\frac{\Delta\bra n,A_{\mu,\nu}n\ket}{\sqrt{\Delta^2\det A_{\mu,\nu}}}
=\frac{|1-\overline n^2\mu|^2-\nu^2}{\sqrt{(1-(|\mu|-\nu)^2)(1-(|\mu|+\nu)^2)}}.
\]
Finally, recalling the definition of $\Delta$, we derive
\[
\det A_{\mu,\nu}=\frac{(1+|\mu|-\nu)(1-|\mu|-\nu)}{(1+|\mu|+\nu)(1-|\mu|+\nu)}
=\frac{(1-\nu)^2-|\mu|^2}{(1+\nu)^2-|\mu|^2}.
\]
Inserting the expressions above into \eqref{gamma}, we obtain \eqref{betaestimate}.
\end{proof}
We now turn to the proof of Theorem~\ref{thm:sharp}.
We let $\tmuo,\tnuo:\R\to\R$ be the bounded, 
piecewise constant, $2\pi$-periodic functions defined 
in $[0,2\pi)$ by
\begin{equation}
\label{tmuo}
\tmuo(\theta)=\begin{cases}
0,&\mathrm{if\ }\theta\in I_1\cup I_3\\
(M-M^{1-2\tau})/(1+M+M^{1-2\tau}+M^{2(1-\tau)}),&\mathrm{if\ }\theta\in I_2\cup I_4
\end{cases}
\end{equation}
and
\begin{equation}
\label{tnuo}
\tnuo(\theta)=\begin{cases}
0,&\mathrm{if\ }\theta\in I_1\cup I_3\\
(M^{2(1-\tau)}-1)/(1+M+M^{1-2\tau}+M^{2(1-\tau)}),&\mathrm{if\ }\theta\in I_2\cup I_4.
\end{cases}
\end{equation}
and we set 
\begin{equation}
\label{tmutnu}
\tmu(z)=-\tmuo(\arg z)\,z\overline z^{-1},\qquad \tnu(z)=-\tnuo(\arg z).
\end{equation}
The following holds.
\begin{Proposition}
\label{prop:sharp}
Let $B$ the unit disk in $\R^2$ and let $\tf\in W^{1,2}(B,\C)$ be defined 
in $B\setminus\{0\}$ by
\[
\tf(z)=|z|^{\iso/c}\left(\Theta_{\tau,1}(\arg z)+i\Theta_{\tau,2}(\arg z)\right).
\]
Then $\tf$ satisfies \eqref{beltrami} with $\mu=\tmu$
and $\nu=\tnu$.
Furthermore, there exists $M_0>1$ such that
\begin{align*}
\beta(\tmu,\tnu)=\frac{\iso}{c},
\end{align*}
for every $M\in(1,M_0^{1/\tau})$ if $\tau>0$
and with no restriction on $M$ if $\tau=0$.
\end{Proposition}
In order to prove Proposition~\ref{prop:sharp}, we first need a lemma.
\begin{Lemma}
\label{lem:mutoA}
Suppose $\mu,\nu$ are of the form~\eqref{mu}--\eqref{nu}
and let $\kf,\ks$ be the corresponding functions defined in \eqref{k}.
Then $A_{\mu,\nu}$ as defined in \eqref{Amu} is given by
\begin{align*}
A_{\mu,\nu}\,(z)=&J(\arg z)\left[
\begin{matrix}\kf(\arg z)&&0\\
0&&\ks(\arg z)
\end{matrix}\right]J^*(\arg z)\\
=&\left[
\begin{matrix}\kf\cos^2\theta+\ks\sin^2\theta &&(\kf-\ks)\sin\theta\cos\theta\\
(\kf-\ks)\sin\theta\cos\theta&&\kf\sin^2\theta+\ks\cos^2\theta
\end{matrix}\right]\\
=&(\kf-\ks)\frac{z\otimes z}{|z|^2}+\ks\mathbf I.
\end{align*} 
\end{Lemma}
\begin{proof}
The assumptions~\eqref{mu}--\eqref{nu} on $\mu,\nu$ imply that
\[
\Delta(z)=(1+\mo(\theta)-\no(\theta))(1-\mo(\theta)-\no(\theta)).
\]
and
\[
\mu(z)=-\mo(\theta)\,(\cos2\theta+i\sin2\theta).
\]
Hence, 
\begin{align*}
\Delta (A_{\mu,\nu})_{11}=&|1-\mu|^2-\nu^2=1+2\mo\cos2\theta+\mo^2-\no^2\\
=&[(1+\mo)^2-\no^2]\cos^2\theta+[(1-\mo)^2-\no^2]\sin^2\theta\\
\Delta (A_{\mu,\nu})_{22}=&|1+\mu|^2-\nu^2\\
=&[(1-\mo)^2-\no^2]\cos^2\theta+[(1+\mo)^2-\no^2]\sin^2\theta\\
\Delta (A_{\mu,\nu})_{12}=&-2\Im(\mu)\\
=&4\mo\sin\theta\cos\theta.
\end{align*}
Dividing by $\Delta$ and observing that
\begin{align*}
\frac{(1+\mo)^2-\no^2}{\Delta}=&\frac{1+\mo+\no}{1-\mo-\no}=\kf\\
\frac{(1-\mo)^2-\no^2}{\Delta}=&\frac{1-\mo+\no}{1+\mo-\no}=\ks\\
\frac{4\mo}{\Delta}=&\kf-\ks,
\end{align*}
we obtain the asserted expression for $A_{\mu,\nu}$.
\end{proof}
\begin{proof}[Proof of Proposition~\ref{prop:sharp}]
By direct check, 
$(\Theta_{\tau,1},\Theta_{\tau,2})$ satisfies \eqref{system} with
$\kf=\tkf$, $\ks=\tks$ and $\alpha_\tau=d/c$.
Hence, in view of Proposition~\ref{prop:polar}, $\tf$ satisfies 
\eqref{beltrami} with $\mu=\tmu$ and $\nu=\tnu$.
In view of Lemma~\ref{lem:reduction} and Lemma~\ref{lem:mutoA},
$\Re(\tf)$ satisfies equation~\eqref{elliptic} with $A=A_\tau$
defined in \eqref{tildeA} and $\Im(\tf)$ satisfies 
equation~\eqref{elliptic} with $A=\widehat{A_\tau}$.
By Theorem~\ref{thm:sharp}--(ii), 
$\Re(\tf)$ and $\Im(\tf)$ are H\"older continuous with exponent exactly
$\beta(\tmu,\tnu)=\gamma(A_\tau)=\gamma(\widehat{A_\tau})$
whenever $M\in(0,M_0^{1/\tau})$ if $\tau>0$ and with no restriction
on $M$ if $\tau=0$.
Thus, Proposition~\ref{prop:sharp} is established.
\end{proof}
\begin{proof}[Proof of Theorem~\ref{thm:sharp}]
The proof is a direct consequence of Proposition~\ref{prop:sharp}.
\end{proof}
\section{Appendix: Reduction to divergence form elliptic equations}
\label{sec:appendix}
We prove the following equivalence result, which implies
Lemma~\ref{lem:reduction} when $\Im(\nu)=0$.
See also \cite{AIM}.
\begin{Lemma}
\label{lem:beltramitoelliptic}
Let $g\in W_{\mathrm{loc}}^{1,2}(\Omega,\C)$ satisfy the Beltrami equation
\begin{equation}
\label{fullbeltrami}
\dbar g=\mu\de g+\nu\overline{\de g}
\qquad\mathrm{in\ }\Omega,
\end{equation}
where $\mu,\nu\in L^{\infty}(\Omega,\C)$ satisfy $|\mu|+|\nu|\le\kappa<1$
a.e.\ in $\Omega$.
Let
$B_{\mu,\nu}$ be the bounded matrix-valued function 
defined in terms of the Beltrami coefficients $\mu,\nu$
by
\begin{align}
\label{Bmunu}
B_{\mu,\nu}=\frac{1}{\Delta_1}\left(\left[\begin{matrix}
|1-\mu|^2&-2\Im(\mu-\nu)\\
-2\Im(\mu+\nu)&|1+\mu|^2
\end{matrix}\right]-|\nu|^2\mathbf{I}\right),
\end{align}
where 
\begin{equation}
\label{Delta1}
\Delta_1=|1+\nu|^2-|\mu|^2
\end{equation}
and let $\widetilde B_{\mu,\nu}$ be defined by
\begin{align}
\label{hatBmunu}
\widetilde B_{\mu,\nu}=\frac{1}{\Delta_2}\left(\left[\begin{matrix}
|1-\mu|^2&-2\Im(\mu+\nu)\\
-2\Im(\mu-\nu)&|1+\mu|^2
\end{matrix}\right]-|\nu|^2\mathbf{I}\right),
\end{align}
where $\Delta_2=|1-\nu|^2-|\mu|^2$.
Then $\Re(g)$ is a weak solution for the elliptic equation~\eqref{elliptic} 
with $A=B_{\mu,\nu}$ and $\Im(g)$ is a weak solution for \eqref{elliptic} with $A=\widetilde B_{\mu,\nu}$.
\end{Lemma}
\begin{proof}
Setting $z=x+iy=(x,y)^T$, $g(z)=u(x,y)+iv(x,y)$,
we have:
\begin{align*}
&\dbar g=\frac{1}{2}\left[\begin{matrix}u_x-v_y\\u_y+v_x\end{matrix}\right]
&&\de g=\frac{1}{2}\left[\begin{matrix}u_x+v_y\\-u_y+v_x\end{matrix}\right].
\end{align*}
Setting
\begin{align*}
&Q=\left[\begin{matrix}0&-1\\1&0\end{matrix}\right]
&&R=\left[\begin{matrix}1&0\\0&-1\end{matrix}\right],
\end{align*}
for every $z$ we have
\begin{align*}
&Qz=\left[\begin{matrix}-y\\x\end{matrix}\right]=iz,
&&Rz=\left[\begin{matrix}x\\-y\end{matrix}\right]=\overline z.
\end{align*}
Hence, we can write
\begin{align*}
&\dbar g=\frac{1}{2}\left(\nabla u+Q\nabla v\right),
&&\de g=\frac{1}{2}R\left(\nabla u-Q\nabla v\right).
\end{align*}
Setting
\begin{align*}
M=\left[\begin{matrix}\Re(\mu)&-\Im(\mu)\\\Im(\mu)&\Re(\mu)\end{matrix}\right],
&&N=\left[\begin{matrix}\Re(\nu)&-\Im(\nu)\\\Im(\nu)&\Re(\nu)\end{matrix}\right],
\end{align*}
equation~\eqref{beltrami} can be written in the form:
\begin{equation*}
\nabla u+Q\nabla v=MR\left(\nabla u-Q\nabla v\right)
+N\left(\nabla u-Q\nabla v\right).
\end{equation*}
It follows that
\[
\left(I-MR-N\right)\nabla u=-\left(I+MR+N\right)Q\nabla v
\]
and consequently $u$ satisfies
\[
\left(I+MR+N\right)^{-1}\left(I-MR-N\right)\nabla u=-Q\nabla v
\]
and $v$ satisfies
\[
-Q\left(I-MR-N\right)^{-1}\left(I+MR+N\right)Q\nabla v=Q\nabla u.
\]
By direct computation,
\begin{align*}
B_{\mu,\nu}=&\left(I+MR+N\right)^{-1}\left(I-MR-N\right)\\
\widetilde B_{\mu,\nu}=&-Q\left(I-MR-N\right)^{-1}\left(I+MR+N\right)Q
=-QB_{-\mu,-\nu}Q.
\end{align*}
Now the conclusion follows observing that $\mathrm{div}(Q\nabla\cdot)=0$.
\end{proof}
We note that the Beltrami coefficients $\mu,\nu$ are uniquely determined by the matrix
$B_{\mu,\nu}=(b_{ij})_{i,j=1,2}$:
\begin{align}
\label{munu}
&\mu=-\frac{b_{11}-b_{22}+i(b_{12}+b_{21})}{1+\mathrm{tr}\,B_{\mu,\nu}+\det B_{\mu,\nu}}
&&\nu=\frac{1-\det B_{\mu,\nu}+i(b_{12}-b_{21})}{1+\mathrm{tr}\,B_{\mu,\nu}+\det B_{\mu,\nu}}.
\end{align}
The formulae above may be obtained as follows. For simplicity, in what follows we denote $B_{\mu,\nu}=B$.
From the definition of $B$ we readily obtain:
\begin{align}
\label{mufinal}
&\Re(\mu)=-\frac{\Delta_1}{4}(b_{11}-b_{22})
&&\Im(\mu)=-\frac{\Delta_1}{4}(b_{12}+b_{21})\\
\label{nu2}
&\Im(\nu)=\frac{\Delta_1}{4}(b_{12}-b_{21})\\
\label{mod}
&1+|\mu|^2-|\nu|^2=\frac{\Delta_1}{2}\mathrm{tr}\,B
\end{align}
The relations~\eqref{mufinal} imply
\begin{equation}
\label{modmu}
|\mu|^2=\frac{\Delta_1^2}{16}\left(\mathrm{tr}\,B^TB-2\det B\right).
\end{equation}
From the definition~\eqref{Delta1} of $\Delta_1$ and \eqref{mod} we derive
\begin{equation}
\label{nu1}
\Re(\nu)=\frac{\Delta_1}{2}\left(1+\frac{\mathrm{tr}\,B}{2}\right)-1.
\end{equation}
From equations~\eqref{nu2}, \eqref{modmu} and \eqref{nu1} we derive
\begin{equation}
\label{modnu}
|\nu|^2=\frac{\Delta_1^2}{16}\left[4(1+\mathrm{tr}\,B)+\mathrm{tr}\,B^TB+2\det B\right]
+1-\frac{\Delta_1}{2}(2+\mathrm{tr}\,B).
\end{equation}
Inserting \eqref{modmu} and \eqref{modnu} into \eqref{mod} we obtain
\[
\Delta_1=\frac{4}{1+\mathrm{tr}\,B+\det B}.
\]
Inserting the expression of $\Delta_1$ above into \eqref{mufinal}, \eqref{nu2} and \eqref{nu1}
we derive the asserted expression \eqref{munu}.  
\begin{proof}[Proof of Lemma~\ref{lem:reduction}]
We need only check that when $\Im(\nu)=0$ we have
\begin{equation}
\label{tildehat}
\widetilde B_{\mu,\nu}=\frac{B_{\mu,\nu}}{\det B_{\mu,\nu}}=\widehat B_{\mu,\nu}.
\end{equation}
Let 
\[
\Gamma_{\mu,\nu}=\left[\begin{matrix}
|1-\mu|^2-\nu^2&-2\Im(\mu)\\
-2\Im(\mu)&|1+\mu|^2-\nu^2
\end{matrix}
\right].
\]
Then
\begin{align*}
B_{\mu,\nu}=\frac{\Gamma_{\mu,\nu}}{\Delta_1},
&&\widetilde B_{\mu,\nu}=\frac{\Gamma_{\mu,\nu}}{\Delta_2}
\end{align*}
with $\Delta_1=(1+\nu)^2-|\mu|^2=(1+\nu+|\mu|)(1+\nu-|\mu|)$
and $\Delta_2=(1-\nu)^2-|\mu|^2=(1-\nu+|\mu|)(1-\nu-|\mu|)$.
On the other hand,
\begin{align*}
&\det\Gamma_{\mu,\nu}=(1+|\mu|+\nu)(1+|\mu|-\nu)(1-|\mu|+\nu)(1-|\mu|-\nu)
\end{align*}
and therefore $\Delta_2=\det\Gamma_{\mu,\nu}/\Delta_1$.
It follows that
\[
\widetilde B_{\mu,\nu}=\frac{\Gamma_{\mu,\nu}}{\Delta_2}
=\frac{\Delta_1}{\det\Gamma_{\mu,\nu}}\Gamma_{\mu,\nu}
=\frac{\Delta_1^2}{\det\Gamma_{\mu,\nu}}\frac{\Gamma_{\mu,\nu}}{\Delta_1}
=\frac{B_{\mu,\nu}}{\det B_{\mu,\nu}},
\]
and \eqref{tildehat} is established.
\end{proof}

\end{document}